\documentstyle{article}

\newtheorem{th}{Theorem}
\newtheorem{de}{Definition}
\newtheorem{co}{Corollary}
\newtheorem{prop}{Proposition}

\def\be{\begin{eqnarray*}}
\def\ee{\end{eqnarray*}}
\def\ben{\begin{eqnarray}}
\def\een{\end{eqnarray}}

\def\R{\mbox{\bf R}}
\def\T{\mbox{\bf T}}
\def\C{\mbox{\bf C}}
\def\Z{{\mbox{\bf Z}}}
\def\Q{{\mbox{\bf Q}}}
\def\CP{{\mbox{\bf CP}}}
\def\Lin{{\rm Lin}_{\bf Z}}
\def\wq{\wedge_Q}
\def\supp{{\rm supp}}
\def\sp{{\rm sp}}

\def\l{\lambda}
\def\Subset{\subset\subset}
\def\O{\Omega}
\def\oo{\omega}

\begin{document}

\title{\huge\bf Almost periodicity in complex analysis
\footnote {Supported by INTAS-99-00089 project} }

\author{\huge Favorov S., Rashkovskii A.}

\date{}

\maketitle

\begin{abstract}
This is a brief survey of up-to-date results on holomorphic almost
periodic functions and mappings in one and several complex
variables, mainly due to the Kharkov mathematical school.

\par

\medskip\noindent


\end{abstract}

\bigskip
While the notion of almost periodic function on $\R$ (or $\R^m$)
seems to be quite understood, it is not the case for {\it
holomorphic} almost periodic functions on a strip in $\C$ (or,
more generally, on a tube domain in $\C^m$). Investigation of the
zero sets of holomorphic almost periodic functions leads to such
objects as almost periodic divisors and holomorphic chains. The
central problem is to determine if a divisor (holomorphic chain)
is generated by a holomorphic almost periodic function (mapping).
Partial results in the one-dimensional situation (\cite{JT},
\cite{KL}, \cite{T1}, \cite{T2}) indicate that the problem is
highly non-trivial. We present here a brief survey of recent
results in this direction, together with some related topics,
mainly due to the Kharkov mathematical school.

\bigskip

{\bf Definitions and basic properties.} We start with standard
notions of the theory of almost periodic functions.

\begin{de}
 A continuous mapping $f$ from $\R^m$ to a metric space $X$ is called almost periodic
if its orbit $\{{\cal T}_tf\}_{t\in\R^m}$ is a relatively compact
set in $C(\R^m,X)$ with respect to the topology of uniform
convergence on $\R^m$, where ${\cal T}_t$ is the translation by
$t\in\R^m$: $({\cal T}_tf)(x)=f(x+t)$.
\end{de}
\begin{de}
A continuous mapping $f$ from a tube domain
$T_\O:=\T^m+i\O=\{z=x+iy:\,x\in\R^m,\,y\in\O\}$ to a metric space
$X$ is called almost periodic on $T_\O$ if its orbit $\{{\cal
T}_tf\}_{t\in\R^m}$ is a relatively compact subset of $C(T_\O,X)$
with respect to the topology of uniform convergence on each tube
subdomain $T_{\O'},\,\O'\Subset\O$.
\end{de}

We reserve the word "function" for mappings to $\C$. Any almost
periodic function $f$ on $T_\O$ has its {\it mean value}
$$
{\cal
M}[f](y)=\lim_{s\to\infty}(2s)^{-m}\int_{[-s,s]^m}f(x+iy)dx,\quad
y\in\O.
$$

The {\it spectrum} of $f$ is the set
$$
\sp f=\{\l\in\R^m:\,a_\l(y):= {\cal M}[f(x+iy)e^{-i\langle
x,\l\rangle}]\not\equiv 0\},
$$
$\langle x,\l\rangle$ being the scalar product in $\R^m$; this set
is at most countable. Then a Fourier series is assigned to each
function $f$,
$$
f\sim\sum_n a_{\l_n}(y)e^{i\langle x,\l_n\rangle};
$$
in addition, if $f$ is holomorphic on a tube domain, then
$$
f\sim\sum_n a_{\l_n}e^{i\langle z,\l_n\rangle},\quad a_{\l_n}\in\C.
$$

Any almost periodic function on a tube domain is uniformly
approximated by
Bochner-Fejer's sums on every subtube $T_{\O'}$, $\O'\Subset\O$:
\be \sum_n \beta_j(\l_n)a_{\l_n}(y)e^{i\langle x,\l_n\rangle}\ee
with $0\le\beta_j<1$, $ \beta_j\to 1$ as $j\to\infty$.

The converse statement is also true: if a function $f(z)$ is
uniformly approximated by finite sums $\sum a_n(y)e^{i\langle
x,\l_n\rangle}$ (or $\sum a_ne^{i\langle z,\l_n\rangle}$) on every
tube subdomain $T_{\O'},\,\O'\Subset\O$, then $f(z)$ is almost
periodic (holomorphic almost periodic) on $T_\O$.

All this can be extended to more general objects (namely, to
distributions), see \cite{Sch}, \cite{R5}. Here we will use the
following notion.

\begin{de} A complex-valued measure $g$ on $T_\O$ is called almost
periodic if for every continuous function $\varphi$ compactly
supported in $T_\O$, the action of $g$ on ${\cal T}_t\varphi$,
$$
g_\varphi(t):=\int {\cal T}_t\varphi\,dg,
$$
is an almost periodic
function in $t\in\R^m$.
\end{de}
Let us remark that the spectrum of $g$ is the union of the spectra
of $g_\varphi$ over all functions $\varphi$.

Since every holomorphic almost periodic function on a tube domain
$T_\O$ has a unique extension to the convex hull of the tube, the
base $\O$ will be always assumed to be convex.

\bigskip

{\bf Zeros of almost periodic functions.} It was proved in
\cite{R1} that the function $\log|f(z)|$ is almost periodic, in
the sense of distributions, for any holomorphic almost periodic
function $f$ on a tube domain. Therefore, the divisor $d_f$ of
such a function $f$ is almost periodic in the sense of
distributions. This means that the current
\ben
\partial\bar\partial\log|f(z)|
=\sum_{j,k}{\partial^2\log|f(z)|\over
\partial z_j\partial\bar z_k}dz_j\wedge d\bar z_k\nonumber
\een
is almost periodic, i.e., all the terms in this sum (that are always
complex-valued measures) are almost periodic. Note that, on the other
hand, almost periodicity of this current does not imply that the function
$f$ itself is almost periodic.

For $m=1$, the divisor $d_f$ of a holomorphic function $f$ on an
open strip $S=T_{(\alpha,\beta)}$ is a positive discrete measure
on $S$. Its mass at each point $a$ equals the multiplicity $k(a)$
of the zero of the function $f$ at this point. The divisor is
almost periodic whenever the convolution
$$
\sum_{a_n\in\,\supp\, d_f}k(a_n)\chi(a_n+t)
$$
is an almost periodic function in the variable $t\in\bf R$ for
every function $\chi\in C(S)$ with $\supp\,\chi\Subset S$. When
$m>1$, almost periodicity of the divisor $d_f$ of a function $f$
holomorphic in $T_\O$ means that the function
$$
\int_{Z_f}k(a)\chi(a+t)
$$
is almost periodic in $t\in\R^m$ for any form $\chi$ of bidegree
$(m-1,m-1)$ with continuous coefficients compactly supported in
$T_\O$; here $Z_f$ is the zero set of $f$.

\begin{th} \label{2} {\rm (\cite{FRR2}, \cite{FRR3}).}
If the trace measure
$$
\sum_{k=1}^m{\partial^2\log|f(z)|\over\partial z_k\partial\bar z_k}
$$
is almost periodic, then so is the divisor $d_f$.
\end{th}

\bigskip

{\bf Jessen's function.} Since $\log|f|$ is almost periodic (in
the sense of distributions), there exists the mean value ${\cal
M}[\log|f|]$. The function
$$
J_f(y):={\cal M}[\log|f|](y)
$$
is called {\it Jessen's function} of $f$. It is convex in $y\in\O$
and generates the Riesz measure
$$
\mu_J=\theta_m\sum_j\partial^2J_f/\partial y_j^2.
$$

\begin{th}\label{3} {\rm (\cite{R1}).} Let $V_f(\O_0,s)$,
$\O_0\Subset\O$, be the volume of the zero set of $f(z)$ inside
the set $\{z:\,x\in [-s,\,s]^m,\ y\in\O_0\}$. Then
$$
\lim_{s\to\infty}(2s)^{-m}V_f(\O_0,s)=\frac{\theta_{2m}}{\theta_m}{\mu_J(\O_0)}
$$
for every $\O_0\Subset\O$ with
$\mu_J(\partial\O_0)=0$. Furthermore, $\mu_J(\O_0)=0$ (i.e.,
$J_f(y)$ is linear on $\O_0$) if and only if $f(z)\neq 0$ on
$T_{\O_0}$.
\end{th}

For $m=1$, the result was proved in \cite{JT}. Moreover, that
seminal work contains a necessary and sufficient (albeit quite
sophisticated) condition for a function to be Jessen's function of
an almost periodic holomorphic function on a strip. As a
consequence, every convex and piece-wise linear function is
Jessen's; when the spectrum of $f$ is a subset of some free
countable additive subgroup of $\R$, the number of linearity
components for $J_f$ is locally finite.

The first statement is no longer true for $m>1$:

\begin{th} {\rm (\cite{R3}).}
A convex piece-wise linear function $J(y)$ on $\O$ is Jessen's
function of an almost periodic holomorphic function $f$ on $T_\O$
if and only if it has the form
$$
J(y)=\sum_{j=1}^\omega \gamma_j\max\{\langle y,\nu^{(j)}\rangle
-h_j ,0\}
+\gamma_0\langle y,\nu^{(0)}\rangle,
$$
where $\omega\le\infty,\:\gamma_j>0,\: \lambda_j\in\R^m,\:
h_j\in\R$; in this case,
$$
f(z)=\prod_{j=1}^\omega f_j\left(\langle z,\nu^{(j)}\rangle
-ih_j\right),
$$
$f_j$ being entire almost periodic functions on $\C$ with real
zeros.
\end{th}

No complete description of Jessen's functions is known in general
situation. Some necessary condition is given in \cite{F2}, which
implies

\begin{th} {\rm (\cite{F2}; for the case of sections of periodic
divisors, see \cite{R7}).} If the spectrum of $f$ is a subset of
some free countable additive subgroup of $\R^m$, then the number
of connected components of the set $\O'\setminus\supp\mu_{J_f}$
(linearity components) is finite for each $\O'\Subset\O$.
\end{th}

The assertion is false in the class of holomorphic almost periodic
functions with arbitrary spectrum.

\bigskip

{\bf Almost periodic divisors.} There exist almost periodic
divisors which cannot be generated by holomorphic almost periodic
functions (for $m=1$, see \cite{T2}; for $m>1$, see \cite{R6}).
This invokes the problem of characterization for divisors of
almost periodic holomorphic functions. The first steps in this
direction were the following results.

\begin{th} \label{t3} {\rm (\cite{T1}).} A skew-symmetric
$N\times N$-matrix $M_d$ with integer entries can be assigned to
each almost periodic divisor $d$ on a strip with the spectrum in
$\Lin\{\l_1,\dots,\l_N\}$, $\l_1,\dots,\l_N\in\R$, such that $d$
is the divisor of a holomorphic almost periodic function on the
strip if and only if  $M_d=0$.
\end{th}

\begin{th} {\rm (\cite{R4})} A skew-symmetric $N\times N$-matrix
$M_D$ with integer entries can be assigned to each $N$-periodic
(with real periods) divisor $D$ on a tube such that
$D$ is the divisor of a holomorphic $N$-periodic function on the tube
if and only if $M_D=0$.
\end{th}

The both matrices $M_d$ and $M_D$ are generated by means of the
increments of $\log \Phi$ for certain functions $\Phi$ constructed
from $d$ and $D$, respectively. A relation between these two
theorems is given by

\begin{th} {\rm (\cite{R6})} Let $d$ be a section of an
$N$-periodic divisor $D$. Then $d$ is almost periodic; moreover,
$d$ is the divisor of some holomorphic almost periodic function if
and only if $D$ is the divisor of some holomorphic $N$-periodic
function.
\end{th}

In order to give a complete description of divisors for almost
periodic holomorphic functions, we need the following notion.

\begin{de}
Bohr's compact set $K$ is the compactification of $\R^m$ with
respect to the topology generated by the sets
$$
\{x\in\mbox{\bf R}^m:|e^{i\langle x,\mu_k\rangle}- e^{i\langle
x_0,\mu_k\rangle}|<\delta,\, 1\le k\le N\}
$$
over all $x_0\in\R^m$, $\delta>0,\,N<\infty,\,\mu_1,\dots,\mu_N
\in\R^m$.
\end{de}

We will consider functions on the "Bohr's tubes" $K_\O:=K+i\O$, so
we need a definition of a function holomorphic on $K_\O$:

\begin{de}
A continuous function $\hat f(\tau+iy)$ is holomorphic on an open
set $\oo\subset K_\O$ if the function ${\cal T}_x\hat f(\tau+iy)$
is a smooth function of $x,\,y\in\R^m$ for small $x$ and
$$\frac{\partial}{\partial\bar z_k}{\cal T}_x\hat
f(\tau+iy)\mid_{x=0}=0, \quad k=1,\dots,m,
$$
for all $\tau+iy\in\oo$.
\end{de}

It is well known that almost periodic functions on $\R^m$ are just
the restrictions of continuous functions on $K$. The corresponding
results for holomorphic almost periodic functions and almost
periodic divisors in tube domain follow.

\begin{prop}{\rm (\cite{F3}).}
Holomorphic almost periodic functions on $T_\O$ are just the
restictions to $T_\O$ of holomorphic functions on $K_\O$.
\end{prop}

\begin{prop}{\rm (\cite{F3}).}
For any almost periodic divisor $d$ on $T_\O$ and any point
$\tau_0+iy_0\in K_\O$, there exists a function $r(\tau+iy)$,
holomorphic on a neighborhood $\oo$ of $\tau_0+iy_0$, with the
following property: for each $x_1+iy_1\in\oo\cap T_\O$, the
restriction of $d$ to some ball in $T_\O$ with center in
$x_1+iy_1$ coincides with the divisor of the restriction of
$r(x+iy)$ to this ball; the function $r(\tau+iy)$ is unique up to
a holomorphic factor without zeros on $\oo$.
\end{prop}

Note that for the divisor of a holomorphic almost periodic
function $f(z)$ we can take $r(\tau+iy)=\hat f(\tau+iy)$; on the
other hand, if the function $r(\tau+iy)$ from the proposition is
well-defined on the whole set $K\times\O$, then $d$ is the divisor
of the holomorphic almost periodic function $r(x+iy)$. Using the
notion of bundle, we obtain the following result.

\begin{th} \label{t1}{\rm (\cite{F3}).}
A line bundle ${\cal F}_d$ over $K_\O$ corresponds to each almost
periodic divisor $d$ on the tube $T_\O$ such that trivial bundles
correspond just to the divisors of almost periodic holomorphic
functions.
\end{th}

It can be proved that the bundle ${\cal F}_d$ is trivial if and
only if the first Chern class of the restriction of the bundle to
the base $K+iy_0$ is trivial. Thus

\begin{th} \label{t2} {\rm (\cite{F1} for $m=1$, \cite{F3} for $m>1$).}
A cohomology class $c(d)\in H^2(K,\Z)$ is assigned to each almost
periodic divisor $d$ on $T_\O$ such that the trivial class
corresponds just to all the divisors of holomorphic almost
periodic functions on $T_\O$. In addition, $c(d)$ is a
homomorphism of the semigroup of all almost periodic divisors to
$H^2(K,\Z)$.
\end{th}

When considering only divisors with spectrum in
$\Lin\{\l_1,\dots,\l_N\}$ and the vectors $\l_1,\dots,\l_N\in\R^m$
linearly independent over $\Z$, Bohr's compact set $K$ in Theorems
\ref{t1} and \ref{t2} can be replaced by the $N$-dimensional torus
$\T^N$. The group $H^2(\T^N,\Z)$ coincides with the group of all
skew-symmetric $N\times N$-matrix with integer entries, so this
particular case of Theorem \ref{t2} is a multidimensional version
of Theorem \ref{t3}.

Theorem \ref{t2} allows to obtain the following sufficient
conditions.

\begin{th} {\rm (\cite{F1} for $m=1$, \cite{F3} for $m>1$).}
If the restriction of an almost periodic divisor $d$ on a tube
$T_\O$ to some tube $T_{\O'}$, $\O'\subset\O$, is the divisor of a
holomorphic almost periodic function $f$ on $T_{\O'}$, then $d$ is
the divisor of a holomorphic almost periodic function $f_1$ on
$T_\O$.
\end{th}

\begin{co} {\rm (\cite{FRR1} for $m=1$, \cite{F3} for $m>1$).}
If the projection of the support $d$ to $\O$ is not dense in $\O$,
then $d$ is the divisor of a holomorphic almost periodic function
$f$ on $T_\O$.
\end{co}

\begin{th} {\rm (\cite{FRR1} for $m=1$, \cite{F3} for $m>1$).}
If an almost periodic divisor $d$ is invariant with respect to the
map $y\mapsto -y$, then $d$ is the divisor of some holomorphic
almost periodic function.
\end{th}

It is not hard to prove that the spectrum $d_f$ for an almost
periodic holomorphic $f$ is contained in the minimal additive
group $G(\sp f)$ containing $\sp\, f$. The converse
statement is also true:

\begin{th} {\rm (\cite{F1} for $m=1$, \cite{F3} for $m>1$).}
If $d$ is the divisor of a holomorphic almost periodic function
$f$, then $d$ is the divisor of some holomorphic almost periodic
function $f_1$ with spectrum in $G(\sp\, d)$.
\end{th}

There exists a convenient representation for the classes $c(d)$.
Let $g(\zeta)$ be an entire function on $\C$ with simple zeros at
all points with integer coordinates. Given $\l,\mu\in\R^m$, denote
by $\l\wq\mu$ the Chern class of the divisor $d^{\l,\mu}$ of the
function
$$
g^{\l,\mu}(z):=g(\langle z,\l\rangle+i\langle z,\mu\rangle),\
z\in\C^m;
$$
it can be proved that $\wq$ is an inner product over $\Q$.
Furthermore, if $\l,\,\mu$ are linearly independent over $\R$,
then $d^{\l,\mu}$ is periodic, and if $\l,\,\mu$ are linearly
dependent over $\R$ and independent over $\Q$, then $d^{\l,\mu}$
is almost periodic.

\begin{th} {\rm (\cite{F3}).}
\label{4} For any almost periodic divisor $d$ on $T_\O$ there
exists a finite number of the standard divisors $d^{\l_k,
\mu_k},\, k=1,\dots,n$, such that
$$
c(d)=\l_1\wq\mu_1+\dots+\l_n\wq\mu_n.
$$
In addition, the divisor
$$
d+\sum_k d^{\mu_k,\l_k}
$$
is the divisor of a holomorphic almost periodic function on
$T_\O$; when $m>1$, all the divisors $d^{\mu_k,\l_k}$ can be taken
periodic.
\end{th}

\bigskip

{\bf Almost periodic holomorphic mappings into affine space.} A
holomorphic mapping $f$ from a tube domain into $\C^q$, $q>1$, is
almost periodic if all its components $f_j$, $1\le j\le q$, are
holomorphic almost periodic functions. It turns out that the zero
set of such a mapping (more precisely, the corresponding
holomorphic chain) need not be almost periodic \cite{FRR3}. The
nature of this phenomenon is that the codimension of the zero set
can drop at infinity. This leads to consideration of the so-called
{\it regular} almost periodic holomorphic mappings introduced in
\cite{R2}:

\begin{de} \label{8}
A holomorphic almost periodic mapping $f:T_\O\to\C^q$, $q\le m$,
is said to be regular if the dimension of the zero set of each
mapping from the closure of the orbit $\{{\cal T}_tf\}_{t\in\R^m}$
of $f$ is at most $m-q$.
\end{de}
A sufficient condition, in terms of spectrum, for a mapping to be
regular \cite{R2} shows that such mappings are, in a sense,
'generic' almost periodic mappings.

The zero set of $f$ can be represented by means of the
Monge-Amp\`ere current $(i\partial\bar\partial \log|f|)^q$, and an
application of the machinery of Monge-Amp\`ere operators gives

\begin{th} {\rm (\cite{FRR2}, \cite{FRR3}).}\label{5.1}
The zero set of any regular almost periodic
holomorphic mapping is almost periodic.
\end{th}

As a result, the zero set has a density $\mu_f$ which is a
non-negative measure on $\Omega$. On the other hand, since the
function $\log|f|$ is almost periodic (in the sense of
distributions) for every almost periodic holomorphic mapping $f$,
one can also define Jessen's function $J_f={\cal M}[\log|f|]$,
which is convex on $\Omega$ (see \cite{Ra1}, \cite{Ra2}).
Nevertheless, there is no direct analog to Theorem \ref{3}: there
exists a regular mapping with no zeros on  $T_\O$ such that the
real Monge-Amp\`ere measure of $J_f$ is strictly positive
\cite{Ra3}. Besides, Jessen's functions of regular mappings have
the following rigidity: if $J_f$ is linear on some open set
$\oo\subset\O$, then $J_f$ is linear everywhere (see \cite{Ra2},
\cite{Ra3}).

A more satisfactory result can be obtained in terms of Jessen's
functions of the components of the mappings. Namely, for a wide
subclass of regular mappings to $\C^q$, $q\le m$ (namely, for
mappings with {\it independent components}), the support of the
real mixed Monge-Amp\`ere measure of $J_{f_1},\ldots,J_{f_q}$
coincides with the union of the supports of $\mu_{g}$ over all
$g(z)=(f_1(z+s^1),\ldots,f_q(z+s^q))$, $s^k\in\C^m$
(see \cite{Ra4}). A sufficient condition for a mapping
to have independent components (in terms of the spectra of the
components) is given in \cite{Ra4}; the condition means that the
spectra are in general position, so such mappings are still
'generic'.

Note that, similar to Theorem \ref{2}, almost periodicity of any
holomorphic chain $Z$ of pure dimension $q$ is equivalent to
almost periodicity of its trace measure $Z\wedge
(i\partial\bar\partial |z|^2)^{m-q}$ {\rm (\cite{FRR2}, \cite{FRR3}).}

\bigskip

{\bf Almost periodic meromorphic functions and
mappings into projective space.} Let $F=(f_0(z):\dots:f_q(z))$ be
a holomorphic almost periodic mapping of a tube $T_\O\subset\C^m$
to the projective space $\CP^q$ equipped with the Fubini-Study
metric. Note that the functions $f_0(z),\dots,f_q(z)$ are not
uniquely defined, but their zeros ("coordinate divisors") are
well-defined; we will suppose that each function is not
identically zero.

In the particular case $q=1$ this gives the class of meromorphic
almost periodic functions; for $m=q=1$ the class was introduced
in \cite{S}. Note that this class is not
closed with respect to either multiplication or addition.

\begin{th} {\rm (\cite{FP1} for $m=1$,
 \cite{F3} for $m>1$)}. The product of two meromorphic almost
periodic functions is almost periodic if and only if the zeros and
poles of the product are uniformly separated in $T_{\O'}$ for any
$\O'\Subset\O$.
\end{th}

\begin{th} {\rm (\cite{FP1} for $m=q=1$,
 \cite{P} for $m=1$, $q>1$,  \cite{F3} for $m>1$, $q>1$)}. There exists a holomorphic
almost periodic mapping from a tube $T_\O\subset\C^m$ to $\CP^q$
with given coordinate divisors $d_0,\dots,d_q$ if and only if

A) all the divisors $d_0,\dots,d_q$ are almost periodic,

B) all the divisors $d_0,\dots,d_q$ have the same Chern class,

C) for any $\O'\Subset\O$, every ball with center in $T_{\O'}$ and
radius $r(\O')$ intersects at most $q$ supports of the divisors
$d_l,\ l=0,\dots,q$.

The mapping can be represented by holomorphic almost periodic
functions $f_0(z),\dots,f_q(z)$ without common zeros if and only
if all the divisors $d_0,\dots,d_q$ have the trivial Chern class.
\end{th}

\bigskip

{\bf Holomorphic functions with almost periodic modulus.} A
complete description for almost periodic divisors on a strip, in
terms of holomorphic functions, is given in the following theorem.

\begin{th}
{\rm (\cite{FRR1})} A
divisor $d$ on a strip $S$ is almost periodic if and only if there
exists a holomorphic function $f$ such that $d_f=d$ and $|f(z)|$
is an almost periodic function on $S$.
\end{th}

It follows from \cite{R2} that if a function $f$, holomorphic on a
tube $T_\O\subset\C^m$,  has almost periodic modulus, then its
divisor is almost periodic and the matrix
$$
M_f:= \left({\rm Im}\,{\cal M}[\partial^2\log|f|/\partial
z_k\bar\partial z_l]\right)_{l,k=1}^m
$$
is zero (since the measures $\partial^2\log|f|/\partial
z_k\bar\partial z_l$ are almost periodic, all the mean values
exist).

The matrix $M_f$ depends only on the Chern class $c(d)$ of the
divisor $d$. Moreover, if the divisors $d$ and $d_0$ have the same
Chern class and $d$ is the divisor of a holomorphic function $f$
with almost periodic modulus, then $d_0$ is the divisor of some
holomorphic function $f_0$ with almost periodic modulus, too.
Using Theorem \ref{4}, we get the following result.

\begin{th}{\rm (\cite{F3})}
If $d$ is an almost periodic divisor in $T_\O$ and the matrix $M_f$
is zero for some $f$ with $d_f=d$, then $d$ is the divisor of a
holomorphic function $f_1$ whose modulus is almost periodic in
$T_\O$.
\end{th}

\bigskip

{\bf Almost periodicity of slices.} It is well-known that a
function, holomorphic and bounded on a strip $S$, is almost
periodic in $S$ if it is almost periodic on at least one straight
line in $S$. In \cite{U1}, \cite{U2} this result was extended to
holomorphic almost periodic functions in tube domains; besides the
usual uniform metric, the integral Stepanoff's, Weyl's, and
Besicovich's metrics were examined as well.

Such a criterion is not true for meromorphic functions in a strip
(which are bounded as mappings to $\CP$). The following theorem is
valid instead.

\begin{th}{\rm (\cite{FP2})}
Let $F(z)$ be a holomorphic mapping from a tube domain $T_\O$ to a
compact complex manifold, uniformly continuous on every domain
$T_\O'$ with $\O'\Subset\O$. If $F(z)$ is almost periodic on some
hyperplane $\R^m+iy_0\subset T_\O$, then $F(z)$ is almost periodic
on $T_\O$.
\end{th}

\bigskip

{\bf Almost periodic solutions of functional equations.} It is
easy to see that a continuous solution of the equation $w^2=f(t)$
with an almost periodic (nonnegative) function $f(t)$ need not be
almost periodic. Nevertheless, the following theorem is true.

\begin{th}{\rm (\cite{BF} for $m=1$, \cite{B2} for $m>1$).}
Let $w(z)$ be a continuous solution of the equation
$$
a_n(z)w^n+a_{n-1}(z)w^{n-1}+\dots+a_1(z)w+a_0(z)=0
$$
on a tube domain $T_\O\subset\C^m$ with holomorphic almost
periodic functions $a_0(z),\dots,a_n(z)$. Then $w(z)$ is
holomorphic almost periodic, too.
\end{th}

In \cite{B1} this theorem was extended to solutions of $F(z,w)=0$
with $F$ holomorphic, almost periodic in $z$ and such that the
mapping $(F(z,w),\,F'_w(z,w))$ is regular in the sense of
Definition \ref{8}.

\bigskip

{\bf Almost periodic functions with spectrum in a cone.} H.\,Bohr
proved that an almost periodic function $f$ on the real axis has
nonnegative spectrum if and only if the function $f$ extends to
the upper half-plane as a bounded holomorphic function; next, the
spectrum of an almost periodic function is bounded if and only if
this function extends to the complex plane as an entire function
of exponential growth.

The following theorems generalize Bohr's results to functions of
several variables.

\begin{th}\label{7}
{\rm(\cite{FU})}. In order that an
almost periodic function $f$ on $\R^m$ have spectrum in a convex
cone $\Gamma\subset\R^m$, it is necessary and sufficient that $f$
extends to the tube domain $T_{\hat\Gamma}$ as a holomorphic
function, bounded in every tube domain $T_{\Gamma'}$.
\end{th}
Here
$$
\hat\Gamma=\{y\in\R^m: \langle y,\,\mu\rangle\ge 0\quad{\rm for\
all} \ \mu\in\Gamma\}
$$
is the conjugate cone to $\Gamma$,  and $\Gamma'$ is an internal
subcone of $\hat\Gamma$.

\begin{th}\label{6}
{\rm(\cite{FU})}. In order that an
almost periodic function $f$ on $\R^m$ have bounded spectrum, it
is necessary and sufficient that $f$ extends to $\C^m$ as an
entire function of the growth at most $C\exp b|z|$. Moreover,
$$
\sup_{x\in\R^m}\limsup_{r\to\infty}{1\over r}\log|f(x+iry)|=
\sup_{\l\in\sp f}\langle-y,\,\l\rangle.
$$
\end{th}

Theorems \ref{7} and \ref{6} are also valid for integral
Stepanoff's metric instead of the usual uniform one.

\bigskip

Department of Mathematics

Kharkov National University

4, pl. Svobody, Kharkov 61077

Ukraine

\medskip

Tek/Nat

H\o gskolen i Stavanger

PB 8002, 4068 Stavanger

Norway

\end{document}